\documentclass{amsart}

\usepackage{amsmath,amssymb,amsthm,a4wide}
\usepackage{graphicx,tikz}
\newtheorem{theorem}{Theorem}

\theoremstyle{definition}

\theoremstyle{remark}

\DeclareMathOperator{\diam}{diam}

\begin{document}

\title[]{ A Stability Version of the Gauss-Lucas Theorem\\ and applications}

\author[]{Stefan Steinerberger}
\address{Department of Mathematics, Yale University, New Haven, CT 06511, USA}
\email{stefan.steinerberger@yale.edu}

\thanks{This work is supported by the NSF (DMS-1763179) and the Alfred P. Sloan Foundation.}

\begin{abstract}
Let $p:\mathbb{C} \rightarrow \mathbb{C}$ be a polynomial. The Gauss-Lucas theorem states that its critical points, $p'(z) = 0$, are contained in the convex hull of its roots. We prove a stability version whose simplest form is as follows: suppose $p$ has $n+m$ roots where $n$ are inside the unit disk,
$$ \max_{1 \leq i \leq n}{|a_i|} \leq 1, \quad \mbox{and $m$ are outside} \quad \min_{n+1 \leq i \leq n+m}{ |a_i|} \geq d > 1 + \frac{2 m}{n},$$
then $p'$ has $n-1$ roots inside the unit disk and $m$ roots at distance at least $(dn - m)/(n+m) > 1$ from the origin and the involved constants are sharp. We also discuss a pairing result: in the setting above, for $n$ sufficiently large each of the $m$ roots has a critical point at distance $\sim n^{-1}$. 
\end{abstract}

\maketitle

\vspace{-15pt}
\section{Introduction and results}
\subsection{Introduction.} The Gauss-Lucas Theorem, first stated by Gauss \cite{gauss} in 1836 and first proved by Lucas \cite{lucas} in 1879, states that if $p:\mathbb{C} \rightarrow \mathbb{C}$ is a polynomial of degree $n$, then the $n-1$ zeroes of $p'$ lie inside the convex hull of the $n$ zeros of $p$. 
This has been refined in various ways \cite{bray, bruj, bruj2, branko, dimitrov, joyal, kalman, malamud, marden, pawlowski, pereira, rav, trevor, siebeck, schm, specht}.  It was recently established by Totik \cite{totik}
that, for sequences of polynomials $p_n$ with $\deg(p_n) \rightarrow \infty$, that if $n-o(n)$ roots of $p$ lie inside a convex domain $K$, then any fixed neighborhood of $K$ contains
$n-o(n)$ roots of $p'$. We prove a sharp, non-asymptotic result in the same spirit.

\begin{theorem} Let $p:\mathbb{C} \rightarrow \mathbb{C}$ be a polynomial having $n$ roots $a_1, \dots, a_n$ inside the unit disk and $m$
roots $a_{n+1}, \dots, a_{n+m}$ outside.  If the roots outside are bounded away from the disk
$$ \min_{n+1 \leq i \leq n+m}{ |a_i|} \geq d > 1 + \frac{2 m}{n},$$
then $p'$ has $n-1$ roots inside the unit disk and $m$ roots of modulus at least $(dn - m)/(n+m) > 1$.
\end{theorem}
Theorem 1 is sharp: consider $p(x) = (x+1)^n (x-d)^m$ for some real $d>1$. The derivative $p'$ has $n-1$ roots in $-1$, $m-1$ roots in $d$ and one root $r$ in
$$ r = \frac{dn -m}{n+m} \qquad \mbox{which requires} \qquad d> 1 + \frac{2m}{n} \quad \mbox{to be outside the unit disk.}$$

\begin{center}
\begin{figure}[h!]
\begin{tikzpicture}[scale=1.2]
\draw [thick] (0,0) circle (0.7cm);
\filldraw [thick] (0,0) circle (0.05cm);
\filldraw [thick] (0.3,0.5) circle (0.05cm);
\filldraw [thick] (-0.2,-0.3) circle (0.05cm);
\filldraw [thick] (0,-0.5) circle (0.05cm);
\filldraw [thick] (1.8,1) circle (0.05cm);
\filldraw [thick] (-1.5,1) circle (0.05cm);
\filldraw [thick] (-1.5,0) circle (0.05cm);
\filldraw [thick] (-0.5,-0.1) circle (0.05cm);
\filldraw [thick] (-0.25,0.6) circle (0.05cm);
\filldraw [thick] (0.5,0) circle (0.05cm);
\filldraw [thick] (0.4,0.1) circle (0.05cm);
\draw [thick] (6,0) circle (0.7cm);
\filldraw [thick] (5.9,0.2) circle (0.05cm);
\filldraw [thick] (6-0.2,-0.1) circle (0.05cm);
\filldraw [thick] (6-0.1,-0.2) circle (0.05cm);
\filldraw [thick] (6.8,0.8) circle (0.05cm);
\filldraw [thick] (6-1.2,0) circle (0.05cm);
\filldraw [thick] (6-1,0.9) circle (0.05cm);
\filldraw [thick] (6-0.6,-0.1) circle (0.05cm);
\filldraw [thick] (6-0.25,0.6) circle (0.05cm);
\filldraw [thick] (6.5,0) circle (0.05cm);
\filldraw [thick] (6.4,0.1) circle (0.05cm);
\node at (0, -1.1) {roots of $p$};
\node at (6, -1.1) {roots of $p'$};
\end{tikzpicture}
\vspace{0pt}
\caption{Roots outside create critical points outside the disk.}
\end{figure}
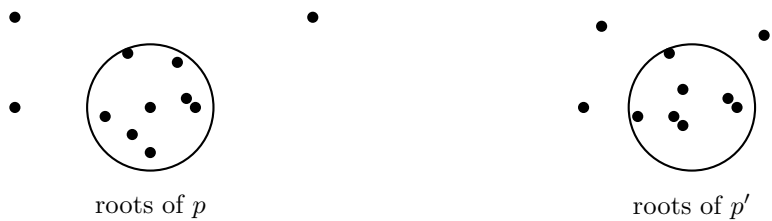
\end{center}

The case $m=1$ is essentially known, though phrased in a somewhat different language, and due to Rahman \cite{rahman} (see also Aziz \cite{aziz} for higher derivatives).
We now formulate a version for general convex domains, $d(z,K)$ denotes the distance from a point $z \in \mathbb{C}$ to a set $K \subset \mathbb{C}$.
\begin{theorem}[Stable Gauss-Lucas] Let $p:\mathbb{C} \rightarrow \mathbb{C}$ be a polynomial having $n$ roots $a_1, \dots, a_n$ inside the convex domain $K \subset \mathbb{C}$ and $m$
roots $a_{n+1}, \dots, a_{n+m}$ outside.  If the roots outside satisfy
$$ \min_{n+1 \leq i \leq n+m}{ d(a_i, K)} \geq 2\diam(K)  \sqrt{\frac{m^2}{n^2} + \frac{m}{n}},$$
then $p'$ has $n-1$ roots close to $K$ in the sense that
$$ d(z,K) \leq \frac{\diam(K) \sqrt{m}}{\sqrt{m+n}}$$
and $m$ other roots satisfying $d(z,K) \geq \diam(K) \sqrt{m/(n+m)}.$
\end{theorem}
We note that the result has a different scaling than Theorem 1 (the scaling in Theorem 1 is a consequence of the boundary of the unit disk having curvature bounded from below). The proof has various degrees of freedom (how to set which parameter etc.) and it is possible to obtain a variety of other results of a similar flavor with the same approach.

\subsection{Pairings of Roots and Critical Points.}
There is recent renewed interest in the interplay between roots of a polynomial and the location of its critical points in the random setting (see e.g. \cite{han, kab, pem}). Kabluchko \cite{kab} has shown that if a random polynomial is constructed by picking its $n$ roots from some probability measure $\mu$, then the roots of the derivative converge to the same measure $\mu$ as $n \rightarrow \infty$. The following recent result is due to O'Rourke \& Williams \cite{or}: suppose $\mu$ is a compactly supported probability measure, $z_1, \dots, z_n$ are randomly drawn from $\mu$ and $\xi \in \mathbb{C}$ is a deterministic point outside the support of $\mu$. If $p$ denotes the random polynomial
$$ p(z) = (z-\xi)\prod_{k=1}^{n}{(z-z_k)},$$
then $p'$ has a root at distance $\sim n^{-1}$ from $\xi$ with high likelihood. We give a deterministic variant.
\begin{center}
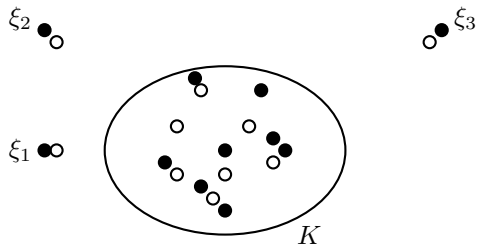
\begin{figure}[h!]
\begin{tikzpicture}[scale = 1.6]
\draw [thick] (0,0) ellipse (1cm and 0.7cm);
\filldraw [thick] (0,0) circle (0.05cm);
\filldraw [thick] (0.3,0.5) circle (0.05cm);
\filldraw [thick] (-0.2,-0.3) circle (0.05cm);
\filldraw [thick] (0,-0.5) circle (0.05cm);
\filldraw [thick] (1.8,1) circle (0.05cm);
\filldraw [thick] (-1.5,1) circle (0.05cm);
\filldraw [thick] (-1.5,0) circle (0.05cm);
\draw [thick] (1.7,0.9) circle (0.05cm);
\draw [thick] (-1.4,0.9) circle (0.05cm);
\draw [thick] (-1.4,0) circle (0.05cm);
\filldraw [thick] (-0.5,-0.1) circle (0.05cm);
\draw [thick] (0.2,0.2) circle (0.05cm);
\draw [thick] (0.4,-0.1) circle (0.05cm);
\draw [thick] (-0.2,0.5) circle (0.05cm);
\draw [thick] (0,-0.2) circle (0.05cm);
\draw [thick] (-0.4,-0.2) circle (0.05cm);
\draw [thick] (-0.4,0.2) circle (0.05cm);
\draw [thick] (-0.1,-0.4) circle (0.05cm);
\node at (0.7,-0.7) {$K$};
\node at (-1.7,0) {$\xi_1$};
\node at (-1.7,1.1) {$\xi_2$};
\node at (2,1.1) {$\xi_3$};
\filldraw [thick] (-0.25,0.6) circle (0.05cm);
\filldraw [thick] (0.5,0) circle (0.05cm);
\filldraw [thick] (0.4,0.1) circle (0.05cm);
\end{tikzpicture}
\caption{Roots outside (filled) have critical points (empty) nearby.}
\end{figure}
\end{center}

\begin{theorem}[Pairing of Roots and Critical Points] Let $K \subset \mathbb{C}$ be a closed convex domain and $\xi_1, \dots, \xi_m \in \mathbb{C} \setminus K$. Then there exists $n_0 \in \mathbb{N}$ and $c>0$, both depending on $m$, $\diam(K)$ and $ \min_{1 \leq i \leq m}{ \min_{z \in K}{\| \xi_i - z\|}} $
such that for all $n> n_0$ and all $z_1, \dots, z_n \in K$, the polynomial
$$ p(z) = \left( \prod_{\ell=1}^{m}(z-\xi_{\ell})\right)\prod_{k=1}^{n}{(z-z_k)}$$
has exactly $m$ critical points outside of $K$ and all of them are at distance $\leq c/n^{}$ from $\left\{\xi_1, \dots, \xi_m\right\}$. Conversely, for each $\xi_i$, there is a critical point in a $c/n^{}-$neighborhood. \end{theorem}

We observe that this is somewhat different from the version of O'Rourke \& Williams \cite{or}: it is completely deterministic but also requires that the $\xi_i$ are not contained inside the convex hull of the support of $\mu$. It turns out that this condition (or some condition in that direction) is necessary in the deterministic setting (as was already observed in \cite{or}): let $K$ to be an annulus centered at the origin containing the boundary of the unit disk, let $\xi=0$ and consider the polynomial
$$ p(z) = (z-\xi)(z^{n}-1) =z \prod_{k=1}^{n}{\left( z - e^{2\pi i k /n}\right)}.$$
For every $\varepsilon > 0$ and and $n$ sufficiently large (depending on $\varepsilon$), all the critical points are contained in the annulus $\left\{z \in \mathbb{C}: 1- \varepsilon \leq |z| \leq 1\right\}$ and none of them are particularly close to $\xi = 0$. The proof of Theorem 3 is fairly explicit and the constant $c$ could be made explicit in terms of everything else if one so desired. The proof also yields the other direction as a byproduct: the critical points either coincide with an $\xi_i$ (which can happen in the case of multiplicity) or are indeed at distance $\gtrsim c_2/n^{}$ away from $\xi_i$ for a constant $c_2$ (depending on the same things as $c$).

\section{Proofs}
\subsection{Proof of Theorem 1}
\begin{proof}
Let $p:\mathbb{C} \rightarrow \mathbb{C}$ be a polynomial of degree $n+m$, having $n$ roots $\left\{a_1, \dots, a_{n+m}\right\}$ where $|a_i| \leq 1$ for all $1 \leq i \leq n$ and 
$$d:= \min\left\{ |a_{n+1}|, |a_{n+2}|, |a_{n+3}|, \dots, |a_{n+m}| \right\} \geq 1.$$
 For our argument it is not important whether the roots are distinct or occur with multiplicity. The derivative $p'$ has $n + m -1$ roots  whose location is determined by the logarithmic derivative
$$ \frac{p'(z)}{p(z)} = \sum_{k=1}^{n+m}{ \frac{1}{ z - a_k}}.$$
We start by obtaining a lower bound on the size of the quantity outside the unit disk caused by the terms contained inside the unit disk. An elementary inequality for real $x > 1$ and complex $|y| \leq 1$ is given by
$$ \mbox{Re} \frac{1}{x - r e^{it}} \geq \frac{1}{x+1} \qquad \mbox{with equality for}~(r,t)=(1,\pi).$$
This implies, for $z = re^{it}$ with $r>1$,
$$ \left|  \sum_{k=1}^{n}{ \frac{1}{ z - a_k}} \right| = \left| \frac{1}{e^{it}} \sum_{k=1}^{n}{ \frac{1}{ z e^{-it} - a_k e^{-it}}} \right|  \geq   \mbox{Re} \sum_{k=1}^{n}{ \frac{1}{ z e^{-it} - a_k e^{-it}}}  \geq \frac{n}{|z| + 1}.$$
We now estimate the size of the remaining term for $|z| < d$. Clearly,
$$ \left|  \sum_{k=n+1}^{n+m}{ \frac{1}{ z - a_k}} \right| \leq  \sum_{k=n+1}^{n+m}{ \frac{1}{ |z - a_k|}} \leq \frac{m}{d - |z|}.$$
In the location of a new root outside the unit disk, the electrostatic forces add up to 0 and we therefore have to have
$$\frac{n}{|z| + 1} \leq  \left|  \sum_{k=1}^{n}{ \frac{1}{ z - a_k}} \right| =  \left|  \sum_{k=n+1}^{n+m}{ \frac{1}{ z - a_k}} \right| \leq  \frac{m}{d - |z|}$$
which then simplifies to
$$ |z| \geq \frac{dn - m}{n+m}.$$
 It remains to show that the derivative still has $m$ roots
outside the unit disk (which, by the preceding argument, are then necessarily at least a controlled distance away from the unit disk). This is done by showing that $p'$ has $n-1$ roots inside the unit disk. We introduce a polynomial collecting all roots inside the unit disk
$$ q(x) = \prod_{k=1}^{n}{(z-a_k)}$$
and observe that critical points of $p$ are zeros of the function
 $$ h(z) = q'(z) + q(z)\sum_{k=n+1}^{n+m}{\frac{1}{z-a_k}}.$$
 $h$ is holomorphic in a neighborhood of the unit disk.
The argument above shows that
$$ |q'(z)| \geq \frac{n}{|z| + 1}|q(z)| \qquad \mbox{for}~z~\mbox{outside the unit disk.}$$
We now apply Rouch\'{e}'s theorem on the boundary of the unit disk (or a slightly larger disk if there happens to a root $|a_i| =1$). Since
$$ d > 1 + \frac{2m}{n} \qquad \mbox{we have} \qquad \frac{n}{|z| + 1} > \frac{m}{|d|-|z|}$$
for $|z| = 1$ and thus
$$ \left|q(z)\sum_{k=n+1}^{n+m}{\frac{1}{z-a_k}} \right|   < | q'(z) |.$$
This shows that the number of roots of $p'$ inside the unit disk, which is the number of roots of $h$ is exactly the same as the number of roots of $q'$ which, by the Gauss-Lucas theorem, is $n-1$.
\end{proof}

\subsection{Proof of Theorem 2}
\begin{proof} The structure of the argument is completely identical to that of Theorem 1, however, some of the computational aspects change.  As before we assume that $a_1, \dots, a_n \in K$ and abbreviate
 $$d:= \min\left\{ d(a_{n+1},K), d(a_{n+2},K), d(a_{n+3},K), \dots, d(a_{n+m},K) \right\}.$$
Let $z \in \mathbb{C} \setminus K$, then we first require a lower bound on
$$ \left|  \sum_{k=1}^{n}{ \frac{1}{ z - a_k}} \right|.$$
By rotational and translational variance, we can again assume that $z \in \mathbb{R}$, that $(0,0) \in K$ is the closest point in $K$ and 
$$z> \sup_{y \in K}{\mbox{Re}~k} = 0.$$ 
A simple computation shows that for all $k \in K$
\begin{align*}
\mbox{Re}\frac{1}{z-k} &= \mbox{Re}\frac{1}{z-k_1 - i k_2} =  \frac{z-k_1}{(z-k_1)^2 + k_2^2}.
\end{align*}
We know that $k_1 \leq 0$ and that $k_1^2 + k_2^2 \leq \diam(K)^2$. An explicit optimization yields that
\begin{align*}
 \frac{z-k_1}{(z-k_1)^2 + k_2^2} &\geq \begin{cases} \frac{d(z,K)}{d(z,K)^2 + \diam(K)^2} \qquad &\mbox{if}~\diam(K) \geq  d(z,K) \\
 \frac{1}{d(z,K) + \diam(K)} \qquad &\mbox{if}~\diam(K) \leq d(z,K)
 \end{cases}\\
 &\geq \frac{d(z,K)}{d(z,K)^2 + \diam(K)^2}.
\end{align*}
Altogether, this implies
$$ \left|  \sum_{k=1}^{n}{ \frac{1}{ z - a_k}} \right| \geq \frac{n \cdot d(z,K)}{d(z,K)^2 + \diam(K)^2} 
\qquad \mbox{for}~z \in \mathbb{C} \setminus K.$$
The upper bound on the effect from outside $K$ is unchanged
$$ \left|  \sum_{k=n+1}^{n+m}{ \frac{1}{ z - a_k}} \right| \leq  \sum_{k=n+1}^{n+m}{ \frac{1}{ |z - a_k|}} \leq \frac{m}{d-d(z,K)}.$$
This shows that the argument in the proof of Theorem 1 is applicable as soon as
$$\frac{n \cdot d(z,K)}{d(z,K)^2 + \diam(K)^2}  >  \frac{m}{d-d(z,K)}$$
This is satisfied as soon as
$$ d > d(z,K)\left( 1 + \frac{m}{n}\right) +   \frac{ \diam(K)^2}{d(z,K)} \frac{m}{n}.$$
Minimizing this quantity shows in the variable $d(z,K)$ shows that we would like to apply it for
$$ d(z,K) = \frac{\diam(K) \sqrt{m}}{\sqrt{m+n}}$$
which ends up requiring
$$ d > 2 \diam(K)  \sqrt{  \frac{m^2}{n^2} + \frac{m}{n}}.$$
\end{proof}

\subsection{Proof of Theorem 3}
\begin{proof} Theorem 2 immediately implies that for $\left\{\xi_1, \dots, \xi_m\right\}$ and $n$ sufficiently large, there are exactly $m$ critical points outside of $K$. It remains
to understand their location. Denoting the roots inside $K$ by $z_1, \dots, z_n$, we obtain that any critical point $z$ satisfies
$$ \left| \sum_{\ell=1}^{n}{\frac{1}{z - z_{\ell}}} \right| = \left| \sum_{\ell=1}^{m}{\frac{1}{\xi_{\ell} - z}} \right| \leq \frac{m}{\min_{1 \leq \ell \leq k}{|z-\xi_{\ell}|}}.$$ 
We know from the proof of Theorem 2 that the left-hand side grows like 
$$ \left| \sum_{\ell=1}^{n}{\frac{1}{z - z_{\ell}}} \right|  \geq \frac{n \cdot d(z,K)}{d(z,K)^2 + \mbox{diam}(K)^2} \qquad \mbox{for}~z \in \mathbb{C} \setminus K.$$
 This requires the right-hand side to be in a $\sim n^{-1}$ neighhborhood of $\left\{ \xi_{1}, \dots, \xi_{m} \right\}$. It also shows that the constant $c$ depend on $m$, the distance of $\left\{\xi_1, \dots, \xi_m\right\}$ to $K$ and the diameter of $K$. It remains to show that a $c/n-$neighborhood of $\xi_1$ necessarily contains a critical point. Critical points satisfy
 $$ \frac{\mbox{mult}(\xi_1)}{z - \xi_1} = \sum_{k=2}^{m}{ \frac{1}{ \xi_{k}-z}} +  \sum_{k=1}^{n}{ \frac{1}{ z_{k}-z}},$$
 where $\mbox{mult}(\xi_1)$ denotes the number of times $\xi_1$ appears in the list and thus w.l.o.g. $\xi_1 \neq \xi_i$ for all $2 \leq i \leq m$. This equation can be rewritten as
 $$ (z - \xi_1)  - \mbox{mult}(\xi_1) \left(  \sum_{k=2}^{m}{ \frac{1}{ \xi_{k}-z}} +  \sum_{k=1}^{n}{ \frac{1}{ z_{k}-z}} \right)^{-1} = 0. \qquad \quad (\diamond)$$
 By construction, the first term in parentheses is bounded by a constant in a sufficiently small neighborhood of $\xi_1$, i.e. there exists $\varepsilon > 0$ such that for all $|z-\xi_1| \leq \varepsilon$ we have
 $$ \left|\sum_{k=2}^{m}{ \frac{1}{ \xi_{k}-z}} \right| \leq C.$$
 The second term is growing linearly in $n$. For $n$ sufficiently large
 $$ \frac{c_1}{n} \leq \left| \mbox{mult}(\xi_1) \left(  \sum_{k=2}^{m}{ \frac{1}{ \xi_{k}-z}} +  \sum_{k=1}^{n}{ \frac{1}{ z_{k}-z}} \right)^{-1} \right| \leq \frac{c_2}{n}$$
 in a sufficiently small neighborhood of $\xi_1$. Applying Rouche's theorem to $(\diamond)$ on the boundary of the disk $\left\{z : |z-\xi_1| = 2 c_2 n^{-1} \right\}$ shows the existence of a critical point in a $2c_2 n^{-1}$ neighborhood. Applying it in the other direction on the boundary of the disk $\left\{z : |z-\xi_1| = c_1 n^{-1}/2 \right\}$ shows that the scaling is optimal, i.e. that the critical point is actually at distance $\sim n^{-1}$ from $\xi_1$.
 \end{proof}

\end{document}